\documentclass[envcountsame]{llncs}
\usepackage[latin1]{inputenc}
\usepackage{amssymb}

\title{{\bf  An Effective Property  of $\omega$-Rational Functions}}

\author{Olivier Finkel }
\institute{Institut de Math\'ematiques de Jussieu - Paris Rive Gauche
 \\  CNRS et Universit\'e Paris 7, France. 
\\ \email{Olivier.Finkel@math.univ-paris-diderot.fr}}

\date{}
\begin{document}

\spnewtheorem{Rem}[theorem]{Remark}{\bfseries}{\itshape}
\spnewtheorem{Exa}[theorem]{Example}{\bfseries}{\itshape}

\spnewtheorem{Pro}[theorem]{Proposition}{\bfseries}{\itshape}
\spnewtheorem{Lem}[theorem]{Lemma}{\bfseries}{\itshape}
\spnewtheorem{Cor}[theorem]{Corollary}{\bfseries}{\itshape}
\spnewtheorem{Deff}[theorem]{Definition}{\bfseries}{\itshape}
\spnewtheorem{Not}[theorem]{Notation}{\bfseries}{\itshape}

\newcommand{\fa}{\forall}
\newcommand{\Ga}{\mathsf\Gamma}
\newcommand{\Gas}{\Ga^\star}
\newcommand{\Gao}{\Ga^\omega}

\newcommand{\Si}{\mathsf\Sigma}
\newcommand{\Sis}{\Si^\star}
\newcommand{\Sio}{\Si^\omega}
\newcommand{\ra}{\rightarrow}
\newcommand{\hs}{\hspace{12mm}

\noi}
\newcommand{\lra}{\leftrightarrow}
\newcommand{\la}{language}
\newcommand{\ite}{\item}
\newcommand{\Lp}{L(\varphi)}
\newcommand{\abs}{\{a, b\}^\star}
\newcommand{\abcs}{\{a, b, c \}^\star}
\newcommand{\ol}{ $\omega$-language}
\newcommand{\orl}{ $\omega$-regular language}
\newcommand{\om}{\omega}
\newcommand{\nl}{\newline}
\newcommand{\noi}{\noindent}
\newcommand{\tla}{\twoheadleftarrow}
\newcommand{\de}{deterministic }
\newcommand{\proo}{\noi {\bf Proof.} }
\newcommand {\ep}{\hfill $\square$}
\renewcommand{\thefootnote}{\star{footnote}} 
\newcommand{\trans}[1]{\stackrel{#1}{\rightarrow}}
\newcommand{\dom}{\operatorname{Dom}}
\newcommand{\im}{\operatorname{Im}}

\maketitle 

\begin{abstract} We prove that  $\om$-regular languages accepted by   B\"uchi or Muller automata satisfy  an effective automata-theoretic version of the Baire property. Then we use this result to obtain a new effective property of rational functions over infinite words which are realized by finite state B\"uchi transducers: for each such function $F: \Sio \ra \Gao$, one can construct a deterministic  B\"uchi automaton $\mathcal{A}$ accepting a dense  ${\bf \Pi}^0_2$-subset of  $\Sio$  such that the restriction of $F$ to $L(\mathcal{A})$ is continuous. 
\noi \end{abstract}

\noi {\bf Keywords.} Decision problems; regular languages of infinite words; infinitary rational relations; omega rational functions;  topology; automatic Baire property; points of continuity.

\section{Introduction}

Infinitary rational relations were firstly studied by Gire and Nivat, \cite{Gire-Phd,Gire-Nivat}. 
The $\om$-rational functions over infinite words, whose  graphs   are (functional)  infinitary rational relations accepted by $2$-tape B\"uchi automata, 
have been studied by several authors \cite{ChoffrutGrigorieffG99,BCPS03,Staiger97,PrieurPhd}. 

In this paper we are mainly interested in the question of the continuity of such $\om$-rational functions. Recall that Prieur proved that 
 one can decide whether a given $\om$-rational function is continuous, \cite{Prieur01,Prieur02}. On the other hand, Carton, Finkel and Simonnet proved that one cannot decide 
whether a given
  $\om$-rational function $f$ has at least one point of continuity,  \cite{CFS08}. Notice that this decision problem is actually $\Sigma_1^1$-complete, hence highly undecidable, \cite{3Applications}.
It was also proved in   \cite{CFS08} that 
 one cannot decide whether the continuity set of a  given
  $\om$-rational function $f$   (its set of continuity points) is a regular (respectively, context-free) $\om$-language.  
 Notice that  the situation was shown to be  quite different in the case of {\it synchronous}  functions. It was  proved in \cite{CFS08} that  
if  $f : A^\omega \to B^\omega$ is an $\om$-rational synchronous function, then the
  continuity set $C(f)$ of~$f$ is $\om$-rational. Moreover 
  If $X$ is  an $\om$-rational ${\bf \Pi}_2^0$ subset of~$A^\omega$, then $X$ is the
  continuity set $C(f)$ of some rational synchronous function~$f$ of
  domain~$A^\omega$.
Notice that these previous works on    the continuity of  $\om$-rational functions had shown that decision problems in this area may be decidable or not, 
(while it is well known that most problems about regular languages accepted by finite automata are decidable). 
 
 We establish in this  paper  a new  effective property of rational functions over infinite words. We first prove that  $\om$-regular languages accepted by   B\"uchi or Muller automata satisfy  an effective automata-theoretic version of the Baire property. Then we use this result to obtain a new effective property of rational functions over infinite words which are realized by finite state B\"uchi transducers: for each such function $F: \Sio \ra \Gao$, on can construct a deterministic  B\"uchi automaton $\mathcal{A}$ accepting a dense  ${\bf \Pi}^0_2$-subset of  $\Sio$  such that the restriction of $F$ to this dense set  $L(\mathcal{A})$ is continuous.

 The paper is organized as follows. We recall basic notions on automata and on the  Borel hierarchy   in Section 2.  The automatic Baire property for regular $\om$-languages is proved in Section 3.  
We prove our main new result on $\om$-rational functions in  Section 4. Some 
concluding remarks are given in Section 5.

\section{Recall of basic notions}

  We assume  the reader to be familiar with the theory of formal ($\om$)-languages  
\cite{Thomas90,Staiger97}.
We recall some  usual notations of formal language theory. 

  When $\Si$ is a finite alphabet, a {\it non-empty finite word} over $\Si$ is any 
sequence $x=a_1\ldots a_k$, where $a_i\in\Si$ 
for $i=1,\ldots ,k$ , and  $k$ is an integer $\geq 1$. The length of $x$ is $|x|=k$.  
 The {\it set of finite words} (including the empty word whose length is zero) over $\Si$ is denoted $\Sis$.
 
  The {\it first infinite ordinal} is $\om$.
 An $\om$-{\it word} over $\Si$ is an $\om$ -sequence $a_1 \ldots a_n \ldots$, where for all 
integers $ i\geq 1$, 
$a_i \in\Si$.  When $\sigma$ is an $\om$-word over $\Si$, we write
 $\sigma =\sigma(1)\sigma(2)\ldots \sigma(n) \ldots $,  where for all $i$,~ $\sigma(i)\in \Si$, and 
$\sigma[n] =\sigma(1)\sigma(2)\ldots \sigma(n)$. 

 The usual concatenation product of two finite words $u$ and $v$ is 
denoted $u\cdot v$ and sometimes just $uv$. This product is extended to the product of a 
finite word $u$ and an $\om$-word $v$: the infinite word $u\cdot v$ is then the $\om$-word such that:
~~  $(u\cdot v)(k)=u(k)$  if $k\leq |u|$ , and 
 $(u\cdot v)(k)=v(k-|u|)$  if $k>|u|$.
  
 The {\it set of } $\om$-{\it words} over  the alphabet $\Si$ is denoted by $\Si^\om$.
An  $\om$-{\it language} over an alphabet $\Si$ is a subset of  $\Si^\om$.

\begin{Deff} : A finite state machine (FSM) is a quadruple $\mathcal{M}=(K,\Si,\delta, q_0)$, where $K$ 
is a finite set of states, $\Sigma$ is a finite input alphabet, $q_0 \in K$ is the initial state
and $\delta$ is a mapping from $K \times   \Si$ into $2^K$. A FSM is called deterministic
 (DFSM) iff: 
$\delta : K \times  \Si \ra K$.

A B\"uchi automaton (BA) is a 5-tuple $\mathcal{A}=(K,\Si,\delta, q_0, F)$ where
  $\mathcal{M}=(K,\Si,\delta, q_0)$
is a finite state machine and $F\subseteq K$ is the set of final states.

A Muller automaton (MA) is a 5-tuple $\mathcal{A}=(K,\Si,\delta, q_0,\mathcal{ F})$ where
$\mathcal{M}=(K,\Si,\delta, q_0)$ is a FSM and $\mathcal{F}\subseteq 2^K$ is the collection of 
designated state sets.

A B\"uchi or Muller automaton is said \de if the associated FSM is deterministic.

Let $\sigma =a_1a_2\ldots a_n\ldots$ be an  $\om$-word over $\Si$.

A sequence of states $r=q_1q_2\ldots q_n\ldots$  is called an (infinite) run of $\mathcal{M}=(K,\Si,\delta, q_0)$ on $\sigma$, 
starting in state $p$, iff:
1) $q_1=p$  and 2) for each $i\geq 1$, $q_{i+1} \in \delta( q_i,a_i)$.

In case a run $r$ of $\mathcal{M}$ on $\sigma$ starts in state $q_0$, we call it simply  ``a run of $\mathcal{M}$ 
on $\sigma$".
For every (infinite) run $r=q_1q_2\ldots q_n\ldots $ of $\mathcal{M}$, $\rm{In}(r)$ is the set of
states in $K$ entered by $M$ infinitely many times during run $r$:
$\rm{In}(r)= \{q\in K  \mid  \exists^\infty i\geq 1 ~  q_i=q\} $ is infinite$\}$.

  For $\mathcal{A}=(K,\Si,\delta, q_0, F)$ a BA,
the \ol~ accepted by $\mathcal{A}$ is:

$L(\mathcal{A})= \{  \sigma\in\Si^\om \mid \mbox{ there exists a  run } r
 \mbox{ of }  \mathcal{A} \mbox{ on }  \sigma \mbox{  such that } \rm{In}(r) \cap F \neq\emptyset \}$.

For $\mathcal{A}=(K,\Si,\delta, q_0, F)$ a MA,
the \ol~ accepted by $\mathcal{A}$ is:

$L(\mathcal{A})= \{  \sigma\in\Si^\om \mid \mbox{ there exists a  run } r
 \mbox{ of }  \mathcal{A} \mbox{ on }  \sigma \mbox{  such that } \rm{In}(r) \in \mathcal{F} \}$.

\end{Deff}

\noi By  R. Mc Naughton's Theorem, see \cite{PerrinPin}, the expressive
 power of deterministic MA (DMA) is equal to the expressive power of non deterministic  MA
(NDMA) which is also equal to the expressive power of non deterministic  BA (NDBA).



\begin{theorem}
  For any $\om$-language  $L\subseteq \Sio$, the following conditions are equivalent:
\begin{enumerate}
\ite   There exists a DMA  that accepts $L$.
\ite   There exists a MA  that accepts $L$.
\ite     There exists a BA  that accepts $L$.
\end{enumerate}

\noi An $\om$-language   $L$ satisfying one of the conditions of the above Theorem is called 
 a regular  $\om$-language. The class of regular  $\om$-languages will
 be denoted by $REG_\om$.
\end{theorem}

\noi  Recall that, from a B\"uchi  (respectively,  Muller)  automaton $\mathcal{A}$, one can effectively construct a deterministic  Muller (respectively, non-deterministic B\"uchi)  automaton  $\mathcal{B}$  such that $L(\mathcal{A})=L(\mathcal{B})$. 

\hs A way to study the complexity of $\om$-languages accepted by various automata 
is to study their topological complexity.

\hs We assume the reader to be familiar with basic notions of topology which
may be found in \cite{Moschovakis80,LescowThomas,Kechris94,Staiger97,PerrinPin}.
If $X$ is a finite alphabet containing at least two letters, then the set  $X^\om$ of infinite words over $X$ 
may be equipped with the product topology of the discrete topology on $X$. 
This topology is induced by  a natural metric  which is called the {\it prefix metric} and defined as follows. 
For $u, v \in X^\om$ and 
$u\neq v$ let $\delta(u, v)=2^{-l_{\mathrm{pref}(u,v)}}$ where $l_{\mathrm{pref}(u,v)}$ 
 is the first integer $n$
such that the $(n+1)^{st}$ letter of $u$ is different from the $(n+1)^{st}$ letter of $v$. 
The topological space  $X^\om$ is a Cantor space. 
 The open sets of $X^\om$ are the sets in the form $W\cdot X^\om$, where $W\subseteq X^\star$.
A set $L\subseteq X^\om$ is a closed set iff its complement $X^\om - L$ is an open set.
Closed sets are characterized by the following:

\begin{Pro}
A set $L\subseteq X^\om$ is a closed set of $X^\om$ iff for every $\sigma\in X^\om$, 

$[\fa n\geq 1,  \exists u\in X^\om$  such that $\sigma (1)\ldots \sigma (n).u \in L]$
 implies that $\sigma\in L$.
\end{Pro}

\noi Define now the next classes of the Borel Hierarchy:

\begin{Deff}
The classes ${\bf \Sigma}_n^0$ and ${\bf \Pi}_n^0$ of the Borel Hierarchy
 on the topological space $X^\om$  are defined as follows:
${\bf \Sigma}^0_1$ is the class of open sets of $X^\om$, ${\bf \Pi}^0_1$ is the class of closed sets 
of $X^\om$.
 And for any integer $n\geq 1$:
${\bf \Sigma}^0_{n+1}$   is the class of countable unions 
of ${\bf \Pi}^0_n$-subsets of  $X^\om$, and 
 ${\bf \Pi}^0_{n+1}$ is the class of countable intersections of 
${\bf \Sigma}^0_n$-subsets of $X^\om$.
\end{Deff}

\begin{Rem} The hierarchy defined above is the  hierarchy of Borel sets of finite rank. 
The Borel Hierarchy is also defined for transfinite levels (see \cite{Moschovakis80,Kechris94}) 
 but we shall not need this  in the sequel. 
\end{Rem}



\noi   It turns out that there is a characterization of ${\bf \Pi}^0_2$-subsets of $X^\om$, involving  
 the notion of $W^\delta$ which we now recall. 

\begin{Deff}
For $W\subseteq X^\star$, we set: $W^\delta=\{\sigma\in X^\om \mid  \exists^\infty i \mbox{ such that } \sigma[i]\in W\}.$
\noi ($\sigma \in W^\delta$ iff $\sigma$ has infinitely many prefixes in $W$).
\end{Deff}

\noi Then we can state the following proposition. 

\begin{Pro}
A subset $L$ of $X^\om$ is a ${\bf \Pi}^0_2$-subset of $X^\om$ iff there exists 
a set $W\subseteq X^\star$ such that $L=W^\delta$.
\end{Pro}

\noi 
It is easy to see, using the 
above characterization of ${\bf \Pi}^0_2$-sets, that every $\om$-language accepted by a deterministic B\"uchi automaton is a 
${\bf \Pi}^0_2$-set.
Thus every regular  $\om$-language is a boolean combination of ${\bf \Pi}^0_2$-sets, 
because it is  accepted by a deterministic Muller  automaton and this implies that  it is a boolean combination 
of $\om$-languages accepted by deterministic B\"uchi automata. 

 Landweber studied  the topological properties of regular  $\om$-languages in \cite{Landweber69}.
He  characterized the regular  $\om$-languages  in each of the Borel classes ${\bf \Sigma}^0_1,  {\bf \Pi}^0_1$, 
${\bf \Sigma}^0_2,  {\bf \Pi}^0_2$, 
 and showed that one can decide, for an effectively given regular  $\om$-language
 $L$, whether $L$ is in ${\bf \Sigma}^0_1,  {\bf \Pi}^0_1$,  ${\bf \Sigma}^0_2$,  or ${\bf \Pi}^0_2$.
 In particular, it turned out that a regular  $\om$-language is in the class ${\bf \Pi}^0_2$ iff it 
is accepted by a deterministic B\"uchi automaton.

 Recall that, from a  B\"uchi or Muller automaton $\mathcal{A}$, one can construct some B\"uchi or Muller automata $\mathcal{B}$ and $\mathcal{C}$, such that 
$L(\mathcal{B})$ is equal to the topological closure  of $L(\mathcal{A})$, and $L(\mathcal{C})$ is equal to the topological  interior of $L(\mathcal{A})$, see \cite{Staiger97,PerrinPin}. 

\section{The automatic Baire property}

In this section we are going to prove an automatic version of the result stating that every Borel (and even every analytic) set has the Baire property. 

  We firstly recall some basic definitions about meager sets, see \cite{Kechris94}. In a topological space $\mathcal{X}$, a set $A\subseteq \mathcal{X}$ is said to be {\it nowhere} dense if 
its closure $\bar{A}$ has empty interior, i.e.  \rm{ Int}($\bar{A}$)$=\emptyset$. A set $A\subseteq \mathcal{X}$ is said to be {\it meager} if it is the union of countably many 
nowhere dense sets, or equivalently if it is included in a countable union of closed sets with empty interiors. This means that $A$ is meager if there exist countably many closed sets $A_n$, $n\geq 1$, such that $A \subseteq \bigcup_{n\geq 1} A_n$  where for every integer $n\geq 1$,   \rm{ Int}($A_n$)$=\emptyset$.  A set is {\it comeager}  if its complement is meager, i.e. if it contains the intersection of countably many dense open sets. Notice that  the notion of {\it meager} set is a notion of {\it small} set, while the notion of  {\it comeager} set is a notion of {\it big} set.

Recall that a Baire space is a  topological space $\mathcal{X}$ in which every  intersection of countably many dense open sets is dense, or equivalently in which every countable union of closed sets with empty interiors  has also an empty interior. 
It is well known that every Cantor space $\Sio$ is a Baire space. In the sequel we will consider only Cantor spaces. 

We now recall the notion of Baire property.  For any sets $A, B \subseteq  \Sio$, we denote $A \Delta B$ the symmetric difference of $A$ and $B$, and we write $A =^\star B$ if and only if 
$A \Delta B$ is meager. 

\begin{Deff}
A set $A \subseteq  \Sio$ has the Baire property (BP) if there exists an open set $U \subseteq  \Sio$ such that  $A =^\star U$. 
\end{Deff}

An important result of descriptive set theory is the following result, see  \cite[page 47]{Kechris94}. 

\begin{theorem}
Every Borel set of a Cantor space has the Baire property. 
\end{theorem}

We now consider regular $\om$-languages $L \subseteq \Sio$ for a finite alphabet $\Si$. These languages are Borel and thus have the Baire property. 
We are going to prove the following automatic version of the above theorem. 

\begin{theorem}\label{aut-Baire}
Let $L = L(\mathcal{A})\subseteq \Sio$  be a  regular $\om$-language accepted by a B\"uchi or Muller automaton $\mathcal{A}$. Then one can construct   B\"uchi automata $\mathcal{B}$ and $\mathcal{C}$ such 
that $L(\mathcal{B})\subseteq \Sio$ is open, $L(\mathcal{C})\subseteq \Sio$ is a countable union of closed sets with empty interior, and  $L(\mathcal{A}) \Delta L(\mathcal{B})\subseteq L(\mathcal{C})$, and we shall say that the $\om$-language $L(\mathcal{A})$ has the automatic   Baire property. 
\end{theorem}

\proo We reason by induction on the topological complexity of the   regular $\om$-language $L = L(\mathcal{A})\subseteq \Sio$ accepted by a B\"uchi automaton  $\mathcal{A}$. 

If  $L = L(\mathcal{A})$ is an open set then we immediately see that we get the result with $\mathcal{B}=\mathcal{A}$ and $\mathcal{C}$ is any B\"uchi automaton accepting the empty set. 

If   $L = L(\mathcal{A})$ is a closed set then $ L \setminus {\rm Int}(L)$ is a closed set with empty interior. Moreover it is known that one can construct from the  B\"uchi automaton 
$\mathcal{A}$ another B\"uchi automaton $\mathcal{B}$ accepting ${\rm Int}(L)$, and then also a   B\"uchi automaton  $\mathcal{C}$ accepting  $ L \setminus {\rm Int}(L)$. Then we have 
 $L(\mathcal{A}) \Delta L(\mathcal{B}) = L \setminus {\rm Int}(L) = L(\mathcal{C})$, with $L(\mathcal{B})$ open and $L(\mathcal{C})$ is a closed set with empty interior. 

We now consider the case of a regular $\om$-language $L = L(\mathcal{A})$ which is a ${\bf \Sigma}^0_2$-set.  Recall that then the  $\om$-language $L$ is accepted by a deterministic  finite automaton  $\mathcal{A}=(K,\Si,\delta, q_0, F)$  with co-B\"uchi acceptance condition. Moreover we   may  assume that the 
 automaton  $\mathcal{A}$ is 
 complete. An $\om$-word $x$ is accepted by $\mathcal{A}$ with  co-B\"uchi acceptance condition if the run of the automaton $\mathcal{A}$ on $x$ (which is then unique since the automaton is deterministic) goes only finitely many times through the set of states $F$. Let now $n\geq 1$ and $L_n$ be the $\om$-language of $\om$-words 
$x\in \Sio$ such that the run of the automaton  $\mathcal{A}$ over $x$ goes at most  $n$ times through the set of states $F$.  We have clearly that $L=\bigcup_{n\geq 1} L_n$. Moreover it is easy to see that for every $n\geq 1$ the set $L_n$ is closed. And the interior of $L_n$ is the union of  basic open sets $u.\Sio$, for $u\in \Sio$,  such that 
the automaton $\mathcal{A}$ enters at most $n$ times in states from $F$ during the reading of $u$ and ends the reading of $u$ in a state $q$ such that none state of $F$ is accessible from this state $q$. We have of course that $L_n \setminus {\rm Int}(L_n)$ is a closed set with empty interior and thus $L_n =^\star {\rm Int}(L_n)$. 
Then $\bigcup_{n\geq 1} L_n =^\star \bigcup_{n\geq 1} {\rm Int}(L_n)$   and   we can easily see that $\bigcup_{n\geq 1} {\rm Int}(L_n)$  is accepted by a 
 B\"uchi  automaton  $\mathcal{B}$ which is essentially 
the automaton $\mathcal{A}$ with a modified   B\"uchi acceptance condition expressing ``some state $q$ has been reached from which none state of $F$ is accessible" .  Moreover $L(\mathcal{A}) \Delta L(\mathcal{B})=( \bigcup_{n\geq 1} L_n ) \Delta ( \bigcup_{n\geq 1} {\rm Int}(L_n) ) \subseteq  \bigcup_{n\geq 1} ( L_n  \setminus {\rm Int}(L_n) )$ and the set $\bigcup_{n\geq 1} ( L_n  \setminus {\rm Int}(L_n) )$ is a countable union of closed sets with empty interiors which is easily seen to be accepted by a co-B\"uchi automaton  $\mathcal{C}$ which is essentially the automaton $\mathcal{A}$  where we have deleted  every state $q$ 
from which none state of $F$ is accessible. 

We now look at boolean operations. Assume that $L(\mathcal{A}) \Delta L(\mathcal{B})\subseteq L(\mathcal{C})$, where $\mathcal{A}$, $\mathcal{B}$, and 
$\mathcal{C}$ are B\"uchi automata, $ L(\mathcal{B})$ is open and $L(\mathcal{C})$ is a countable union of closed sets with empty interiors.   Notice that  
$(\Sio \setminus L(\mathcal{A})) \Delta (\Sio \setminus L(\mathcal{B})) =L(\mathcal{A}) \Delta L(\mathcal{B})$. Moreover we can construct a B\"uchi automaton $\mathcal{D}$ accepting the closed set 
$(\Sio \setminus L(\mathcal{B}))$ and next also a B\"uchi automaton $\mathcal{D}'$ accepting the open set   ${\rm Int}((\Sio \setminus L(\mathcal{B}))$, and  
a B\"uchi automaton
$\mathcal{E}$ accepting the closed set with empty interior $( L(\mathcal{D})  \setminus  {\rm Int}(L(\mathcal{D}) )$. 
It is now easy to see that $(\Sio \setminus L(\mathcal{A})) \Delta L(\mathcal{D}')  \subseteq  (\Sio \setminus L(\mathcal{A})) \Delta (\Sio \setminus L(\mathcal{B}))   \cup L(\mathcal{E}) 
 \subseteq  L(\mathcal{C}) \cup L(\mathcal{E})$ and we can construct a  B\"uchi automaton $\mathcal{C}'$ accepting the $\om$-language $ L(\mathcal{C}) \cup L(\mathcal{E})$ which is a 
  countable union of closed sets with empty interiors so that $(\Sio \setminus L(\mathcal{A})) \Delta  L(\mathcal{D}')  \subseteq L(\mathcal{C}')$. Notice that this implies that the automatic Baire property stated in the theorem is satisfied for regular $\om$-languages in the Borel class ${\bf \Pi}^0_2$ since we have already solved the case of the class ${\bf \Sigma}^0_2$. 
  
   We now consider the finite union operation. Assume we have   $L(\mathcal{A}) \Delta L(\mathcal{B})\subseteq L(\mathcal{C})$, and  $L(\mathcal{A}') \Delta L(\mathcal{B}')\subseteq L(\mathcal{C}')$ where $\mathcal{A}$,  $\mathcal{A}'$, $\mathcal{B}$,  $\mathcal{B}'$, and 
$\mathcal{C}$, $\mathcal{C}'$,  are B\"uchi automata, $ L(\mathcal{B})$ and $ L(\mathcal{B}')$ are  open and $L(\mathcal{C})$ and $L(\mathcal{C}')$ are  countable unions of closed sets with empty interiors.  Now we can see that $ ( L(\mathcal{A}) \cup  L(\mathcal{A}') )  \Delta  ( L(\mathcal{B}) \cup  L(\mathcal{B}') )  \subseteq (L(\mathcal{A}) \Delta L(\mathcal{B})) 
\cup (L(\mathcal{A}') \Delta L(\mathcal{B}'))   \subseteq  L(\mathcal{C})  \cup L(\mathcal{C}')$.  Moreover we can construct  B\"uchi automata $\mathcal{B}''$ and $\mathcal{C}''$ such that 
$L(\mathcal{B}'')$ is the open set $ L(\mathcal{B}) \cup  L(\mathcal{B}')$ and $L(\mathcal{C}'')=L(\mathcal{C})  \cup L(\mathcal{C}')$ is a   countable union of closed sets with empty interiors and then we have $ ( L(\mathcal{A}) \cup  L(\mathcal{A}') )  \Delta   L(\mathcal{B}'')  \subseteq L(\mathcal{C}'')$. 
      
We now return to the general case of a regular $\om$-language $L \subseteq  \Sio$, accepted by a B\"uchi or Muller automaton. We know that we can construct a deterministic Muller automaton  $\mathcal{A}=(K,\Si,\delta, q_0,\mathcal{ F})$  accepting $L$. Recall that   $\mathcal{F}\subseteq 2^K$ is here the collection of 
designated state sets.  For each state $q\in K$, we now denote by    $\mathcal{A}_q$ the automaton $\mathcal{A}$ but viewed as a (deterministic) B\"uchi automaton with the  single accepting state $q$, i.e.  $\mathcal{A}_q=(K,\Si,\delta, q_0,\{ q \})$. We know that the languages $L(\mathcal{A}_q)$  are Borel 
 ${\bf \Pi}^0_2$-sets and thus satisfy the automatic Baire property. 
Moreover we  have the following equality: 
$$ L(\mathcal{A}) = \bigcup_{F\in \mathcal{ F} } [ \cap_{q\in F} L(\mathcal{A}_q) \setminus \cup_{q\notin F} L(\mathcal{A}_q) ]$$

This implies, from the case of  ${\bf \Pi}^0_2$ $\om$-regular languages and from the preceding remarks about the preservation of the automatic Baire property by boolean operations, that we can construct B\"uchi  automata $\mathcal{B}$ and $\mathcal{C}$, such that $ L(\mathcal{B})$ is open and $L(\mathcal{C})$ is a countable union of closed sets with empty interiors, which satisfy $L(\mathcal{A}) \Delta L(\mathcal{B})\subseteq L(\mathcal{C})$. 
\ep 

\begin{Cor}
On can decide, for a given  B\"uchi  or Muller automaton $\mathcal{A}$, whether $L(\mathcal{A})$ is meager. 
\end{Cor}

\proo Let $\mathcal{A}$ be a B\"uchi  or Muller automaton. The $\om$-language $L(\mathcal{A})$ has the  automatic Baire property and we can construct 
B\"uchi  automata $\mathcal{B}$ and $\mathcal{C}$, such that $ L(\mathcal{B})$ is open and $L(\mathcal{C})$ is a countable union of closed sets with empty interiors, which satisfy $L(\mathcal{A}) \Delta L(\mathcal{B})\subseteq L(\mathcal{C})$. It is easy to see that  $L(\mathcal{A}) $ is meager if and only if 
$L(\mathcal{B})$ is empty, since any non-empty open set is non-meager, and this can be decided from the automaton $\mathcal{B}$.  \ep 

\begin{Rem}
The above Corollary followed already from Staiger's paper \cite{Staiger98}, see also \cite{MMS17}.  So we get here another proof of this result, based on the automatic Baire property. 
\end{Rem}

\section{An application to $\om$-rational functions}

\subsection{Infinitary rational relations}

\noi We now  recall the definition  of infinitary rational relations,  via definition by B\"uchi transducers:

\begin{Deff}
  A $2$-tape B\"uchi automaton is a sextuple $\mathcal{T}=(K, \Si, \Ga, \Delta,
  q_0, F)$, where $K$ is a finite set of states, $\Si$ and $\Ga$ are finite
  sets called the input and the output alphabets, $\Delta$ is a finite
  subset of $K \times \Sis \times \Gas \times K$ called the set of
  transitions, $q_0$ is the initial state, and $F \subseteq K$ is the set
  of accepting states.  \nl A computation $\mathcal{C}$ of the automaton
  $\mathcal{T}$ is an infinite sequence of consecutive transitions
  \begin{displaymath}
  (q_0, u_1, v_1, q_1), (q_1, u_2, v_2, q_2), \ldots 
  (q_{i-1}, u_{i}, v_{i}, q_{i}),  (q_i, u_{i+1}, v_{i+1}, q_{i+1}), \ldots 
  \end{displaymath}
  The computation is said to be successful iff there exists a final state
  $q_f \in F$ and infinitely many integers $i\geq 0$ such that $q_i=q_f$.
  The input word and output word of the computation are respectively
  $u=u_1.u_2.u_3 \ldots$ and $v=v_1.v_2.v_3 \ldots$ The input and the
  output words may be finite or infinite.  The infinitary rational relation
  $R(\mathcal{T})\subseteq \Sio \times \Ga^\om$ accepted by the  $2$-tape B\"uchi automaton  $\mathcal{T}$ is the set of couples $(u, v) \in \Sio \times
  \Ga^\om$ such that $u$ and $v$ are the input and the output words of some
  successful computation $\mathcal{C}$ of $\mathcal{T}$.  The set of
  infinitary rational relations will be denoted $RAT_2$.
\end{Deff} 

\noi If $R(\mathcal{T})\subseteq \Sio \times \Ga^\om$ is an infinitary rational
relation recognized by the  $2$-tape B\"uchi automaton  $\mathcal{T}$ then we denote
$$Dom(R(\mathcal{T}))=\{ u \in \Sio \mid \exists v \in \Gao ~~(u, v) \in
R(\mathcal{T}) \}$$
\noi
 and 
$$Im(R(\mathcal{T}))=\{ v \in \Gao \mid \exists u \in
\Sio (u, v) \in R(\mathcal{T}) \}.$$

\noi It is well known that, for each infinitary rational relation
$R(\mathcal{T})\subseteq \Sio \times \Ga^\om$, the sets
$Dom(R(\mathcal{T}))$ and $Im(R(\mathcal{T}))$ are regular $\om$-languages and that one can construct, from the 
B\"uchi transducer $\mathcal{T}$, some B\"uchi automata $\mathcal{A}$ and $\mathcal{B}$ accepting the $\om$-languages 
$Dom(R(\mathcal{T}))$ and $Im(R(\mathcal{T}))$.

\hs The  $2$-tape B\"uchi automaton  $\mathcal{T}=(K, \Si, \Ga, \Delta, q_0, F)$ is said
to be synchronous if the set of transitions $\Delta$ is a finite subset of
$K \times \Si \times \Ga \times K$, i.e. if each transition is labelled
with a pair $(a, b)\in \Si \times \Ga$.  An infinitary rational relation
recognized by a synchronous  $2$-tape B\"uchi automaton  is in fact an $\om$-language
over the product alphabet $\Si \times \Ga$ which is accepted by a B\"uchi
automaton. It is called a synchronous infinitary rational relation. An
infinitary rational relation is said to be asynchronous if it can not be
recognized by any synchronous  $2$-tape B\"uchi automaton.  Recall now the following
undecidability result of C. Frougny and J. Sakarovitch.

\begin{theorem}[\cite{FrougnySakarovitch93}]
  One cannot decide whether a given infinitary rational relation is
  synchronous.
\end{theorem}

 We proved in \cite{Fin-HI} that many decision problems about infinitary rational relations are highly undecidable. In fact many of them, like 
the universality problem, the equivalence problem, the inclusion problem, the cofiniteness problem, the unambiguity problem,  
are $\Pi_2^1$-complete, hence located at the second level of the analytical hierarchy.

\subsection{Continuity  of   $\om$-rational functions}

Recall that an infinitary rational relation $R(\mathcal{T})\subseteq \Sio \times \Ga^\om$ is said to be functional iff it is the graph of a function,  i.e. iff 
$$[\fa x \in Dom(R(\mathcal{T})) ~~ \exists ! y \in Im(R(\mathcal{T})) ~~ (x, y) \in R(\mathcal{T})].$$ 

\noi Then the  functional relation $R(\mathcal{T})$  defines an $\om$-rational (partial) function
$F_{\mathcal{T}}: Dom(R(\mathcal{T})) \subseteq \Sio \ra \Gao$ by:
 for each $u\in Dom(R(\mathcal{T}))$, $F_{\mathcal{T}}(u)$ is the unique
$v\in\Gao$ such that $(u, v) \in R(\mathcal{T})$. 

\noi An $\om$-rational
(partial) function $f: \Sio \ra \Gao$ is said to be synchronous if there is
a synchronous  $2$-tape B\"uchi automaton  $\mathcal{T}$ such that
$f=f_{\mathcal{T}}$.  \nl An $\om$-rational (partial) function $f: \Sio \ra
\Gao$ is said to be asynchronous if there is no synchronous $2$-tape B\"uchi automaton  $\mathcal{T}$ such that $f=f_{\mathcal{T}}$.

\hs Recall the following previous decidability result. 

\begin{theorem}[\cite{Gire86}]
  One can decide whether an infinitary rational relation recognized by a
  given  $2$-tape B\"uchi automaton  $\mathcal{T}$ is a functional infinitary
  rational relation (respectively, a synchronous functional infinitary
  rational relation).

\end{theorem}


    It is  very natural to consider the notion of continuity for  $\om$-rational  functions defined by  $2$-tape B\"uchi automata. 

 We recall that a function $f: Dom(f) \subseteq
\Sio \ra \Gao$, whose domain is $Dom(f)$, is said to be continuous at point
$x\in Dom(f)$ if :
$$\forall n\geq 1 ~~~\exists k \geq 1 ~~~\forall y\in Dom(f) ~~~[~ d(x, y)<2^{-k} \Rightarrow d(f(x), f(y))<2^{-n}~] $$

\noi  The continuity set $C(f)$ of the function $f$ is the set of points of
continuity of $f$. Notice that the continuity set $C(f)$ of a function $f: \Sio \ra \Gao$ is always a Borel ${\bf \Pi}^0_2$-subset of $\Sio$, see \cite{CFS08}. 

\hs The function $f$ is said to be continuous if it is continuous at every
point $x\in Dom(f)$, i. e. if $C(f)=Dom(f)$.

\hs Prieur proved the following decidability result. 

\begin{theorem}[Prieur \cite{Prieur01,Prieur02}]
   One can decide whether a given $\om$-rational function is continuous. 
\end{theorem}

\noi On the other hand   the following undecidability result was proved in  \cite{CFS08}. 

\begin{theorem}[see  \cite{CFS08}]
  One cannot decide whether a given
  $\om$-rational function $f$ has at least one point of continuity.
\end{theorem}

\noi The exact complexity of this undecidable problem was given in \cite{3Applications}.  It is  $\Sigma_1^1$-complete   to    determine  whether a given    
$\om$-rational function $f$ has at least one point of continuity. 

\hs  We now consider  the continuity set of an $\om$-rational function and its possible complexity. 
The following undecidability result was proved in  \cite{CFS08}. 

\begin{theorem}[see  \cite{CFS08}]\label{reg}
  One cannot decide whether the continuity set of a  given
  $\om$-rational function $f$ is a regular (respectively, context-free) $\om$-language.  
\end{theorem}

\noi The exact complexity of the first above  undecidable problem, and an approximation of the complexity  of the second one,  were given in \cite{3Applications}.  
 It is  $\Pi_1^1$-complete   to    determine  whether the continuity set $C(f)$ of a given    
$\om$-rational function $f$ is a regular $\om$-language. Moreover the problem to determine  whether the continuity set $C(f)$ of a given    
$\om$-rational function $f$  is a context-free $\om$-language  is $\Pi_1^1$-hard 
and in the class $\Pi_2^1 \setminus \Sigma_1^1$.

\hs The situation is quite different in the case of {\it synchronous}  functions. The following results were proved in \cite{CFS08}.

\begin{theorem}[\cite{CFS08}]
  Let $f : A^\omega \to B^\omega$ be a rational synchronous function.  The
  continuity set $C(f)$ of~$f$ is rational.
\end{theorem}

\begin{theorem}[\cite{CFS08}] \label{thm:charContSet}
  Let $X$ be a rational ${\bf \Pi}_2^0$ subset of~$A^\omega$.  Then $X$ is the
  continuity set $C(f)$ of some rational synchronous function~$f$ of
  domain~$A^\omega$.
\end{theorem}

We are now going to  prove another  effective result about $\om$-rational  functions. 

We first recall the following result of descriptive set theory, in the particular case of Cantor spaces $\Sio$ and $\Gao$. A Borel function $f: \Sio \ra \Gao$ is a function for which the inverse image of any Borel subset of $\Gao$, or equivalently of any open set of $\Gao$, is a Borel subset of $\Sio$. 

\begin{theorem}[see Theorem 8.38  of \cite{Kechris94}]\label{cl}
Let $\Si$ and $\Ga$ be two finite alphabets and  $f: \Sio \ra \Gao$ be a Borel function. Then there is a dense ${\bf \Pi}^0_2$-subset $G$ of $\Sio$ such that the restriction of $f$ to $G$ is continuous. 
\end{theorem}

We now state an automatic version of this theorem. 

\begin{theorem}\label{aut}
Let $\Si$ and $\Ga$ be two finite alphabets and  $f: \Sio \ra \Gao$ be an  $\om$-rational   function. Then there is  a dense $\om$-regular ${\bf \Pi}^0_2$-subset $G$ of $\Sio$ such that the restriction of $f$ to $G$ is continuous. 
Moreover one can construct, from a $2$-tape B\"uchi  automaton accepting the graph of the function $f$,  a deterministic  B\"uchi automaton accepting a dense ${\bf \Pi}^0_2$-subset $G$ of $\Sio$ such that the restriction of $f$ to $G$ is continuous. 
\end{theorem}

\proo
In the classical context of descriptive set theory, the proof of the above  Theorem \ref{cl}  is the following.  Let $U_n$ be a basic open subset of $\Gao$ (there is a countable basis for the usual Cantor topology on $\Gao$). Then $f^{-1}(U_n)=R_n$ is a Borel subset 
of $\Sio$ since $f$ is a Borel function. Thus $R_n$ has the Baire property and there exists some open set $V_n$ such that $R_n \Delta V_n \subseteq F_n$, where $F_n$ ia countable union of closed sets with empty interiors. Let now $G_n=\Sio \setminus F_n$ and $G= \bigcap_{n\geq 1} G_n = \Sio \setminus  \bigcup_{n\geq 1} F_n $. Then $G$ is countable intersection of dense open subsets of $\Sio$, hence also a dense ${\bf \Pi}^0_2$-subset $G$ of $\Sio$ since in the Cantor space any  countable intersection of dense open sets is dense. Moreover the restriction $f_G$ of the function $f$ to $G$ is continuous. Indeed  for every basic open set $U_n$ it holds that $f_G^{-1}(U_n) = f^{-1}(U_n)   \cap G = V_n \cap G$ is an open subset of $G$, and this implies that the inverse image of any open subset 
of $\Gao$ by $f_G$ is also an open subset of $G$. 

\hs 

  We now want to reason in the automata-theoretic framework in order to prove  Theorem \ref{aut}. Let then $\Si$ and $\Ga$ be two finite alphabets and  $f: \Sio \ra \Gao$ be an $\om$-rational   function whose graph is accepted by a  $2$-tape B\"uchi  automaton $\mathcal{A}=(K, \Si, \Ga, \Delta, q_0, F)$.

 We assume that we have an enumeration of the finite words over the alphabet 
$\Ga$ given by $(u_n)_{ n\geq 1}$, $u_n\in \Gas$.  For $q\in K$ we also denote $\mathcal{A}_q$ the automaton $\mathcal{A}$ in which we have changed the initial state so that the initial state of  $\mathcal{A}_q$  is $q$ instead of $q_0$. 

Let us  now  consider the basic open set of the space $\Gao$ given by  $U_n=u_n\cdot \Gao$. We first describe $f^{-1}(U_n)$.
An $\om$-word $x\in\Sio$ belongs to the set 
$f^{-1}(U_n)$ iff $x$ can be written in the form $x=v\cdot y$ for some words $v\in \Sis$ and $y\in \Sio$, and there is a partial run of  the automaton  $\mathcal{A}$ reading $(v, u_n)$ for which $\mathcal{A}$  is in state $q$ after having read the initial pair  $(v, u_n) \in \Si^\star \times \Ga^\star$ (where the finite words $v$ and $u_n$ might have different lengths if the automaton $\mathcal{A}$ is not synchronous), and  $y\in Dom(L(\mathcal{A}_q))$. Recall that $L(\mathcal{A}_q)\subseteq  (\Si \times \Ga)^\om$ is an infinitary rational relation and that $ Dom(L(\mathcal{A}_q))$ is then a regular $\om$-language and that one can construct from  $\mathcal{A}$ a deterministic Muller automaton accepting this $\om$-language  $Dom(L(\mathcal{A}_q))$ which will be denoted $L_q$. We also denote $T(u_n, q)$ the set of finite words $v$ over $\Si$  such that the automaton 
 $\mathcal{A}$ is in state $q$ after having read the initial pair $(v, u_n) \in \Si^\star\times \Ga^\star$. Then the following equality   holds: 
 $$f^{-1}(U_n) = \bigcup_{q\in K} T(u_n, q) \cdot L_q$$
\noi We can now apply the automatic Baire property stated in the above Theorem \ref{aut-Baire}. Then for each regular $\om$-language $L_q$, one can construct 
a deterministic Muller automaton accepting an open set $O_q$ and a deterministic Muller automaton accepting a countable union $W_q$ of closed sets with empty interiors, such that 
for each $q\in K$, 
$$L_q \Delta O_q  \subseteq W_q$$
Now we set 
$$V_n =  \bigcup_{q\in K} T(u_n, q) \cdot O_q$$
and 
$$F_n =  \bigcup_{q\in K} T(u_n, q) \cdot W_q$$
\noi  Notice that each set  $T(u_n, q)$  is countable  and that for each finite word $u\in T(u_n, q)$    it is easy to see that the set $u\cdot O_q$ is open and that the set 
$u\cdot  W_q$ is a countable union of closed sets with empty interiors. 
Thus it is easy to see that $V_n$ is open, and that $F_n$ is a countable union  of closed sets with empty interiors. 
Moreover it is easy to see that $V_n$ and $F_n$ are regular $\om$-languages since each set $T(u_n, q)$ is a regular language of finite words over the alphabet $\Si$. 
  Moreover   it holds that:
$$f^{-1}(U_n) \Delta V_n \subseteq F_n$$
We now prove that $F=\bigcup_{n\geq 1}F_n$ is itself a regular $\om$-language. It holds that 
$$ F = \bigcup_{n\geq 1}F_n = \bigcup_{n\geq 1} \bigcup_{q\in K} T(u_n, q) \cdot W_q =  \bigcup_{q\in K} \bigcup_{n\geq 1} T(u_n, q) \cdot W_q$$
\noi Consider now the $2$-tape  automaton $\mathcal{B}_q$ which is like the $2$-tape automaton $\mathcal{A}$ but reads only pairs of finite words in $\Si^\star \times \Ga^\star$ and has the state $q$ as  unique accepting state. Let then $\mathcal{C}_q$ be a finite  automaton which reads only finite words over the alphabet $\Si$ and such that $L(\mathcal{C}_q)={\rm Proj}_{\Sis}(L(\mathcal{B}_q))$ is the projection of the language $L(\mathcal{B}_q)$ on $\Sis$. We can construct, from the automaton $\mathcal{A}$, the  automata  $\mathcal{B}_q$ and 
$\mathcal{C}_q$ for each $q\in K$. Now it holds that: 
$$ F = \bigcup_{n\geq 1}F_n =  \bigcup_{q\in K} \bigcup_{n\geq 1} T(u_n, q) \cdot W_q = \bigcup_{q\in K} L(\mathcal{C}_q)\cdot W_q $$
\noi On the other hand, for each finite word $u\in \Sis$, the set $u \cdot W_q $ is a  countable union of closed sets with empty interiors, since $W_q$ is a countable union of closed sets with empty interiors. Thus the set 
$$ F = \bigcup_{q\in K} L(\mathcal{C}_q)\cdot W_q $$
\noi is also a countable union of closed sets with empty interiors, since $K$ is finite and each language $ L(\mathcal{C}_q)$ is countable. Moreover the $\om$-language $F$ is regular and 
we can construct,  from the automata $\mathcal{C}_q$ and from the deterministic Muller automata accepting the $\om$-languages $W_q$, a deterministic Muller automaton accepting $F$. 

We can now reason as in the classical case by setting $G_n=\Sio \setminus F_n$ and $G= \bigcap_{n\geq 1} G_n = \Sio \setminus  \bigcup_{n\geq 1} F_n  = \Sio \setminus F$. Then $G$ is a countable intersection of dense open subsets of $\Sio$, hence also a dense ${\bf \Pi}^0_2$-subset $G$ of $\Sio$.  Moreover we can construct a deterministic Muller automaton and even a 
deterministic B\"uchi automaton (since $G$ is a ${\bf \Pi}^0_2$-set, see \cite[page 41]{PerrinPin}) accepting $G$, and 
the restriction $f_G$ of the function $f$ to $G$ is continuous.
\ep

\begin{Rem}
The above  dense ${\bf \Pi}^0_2$-subset $G$ of $\Sio$ is {\bf comeager} and thus Theorem \ref{aut} shows that  one can construct a deterministic B\"uchi automaton accepting   a ``{\bf large}'' $\om$-rational subset of $\Sio$ on which the function $f$ is continuous. 
\end{Rem}

\section{Concluding remarks}

We have proved a new effective property of $\om$-rational functions. We hope this property will be useful for further studies involving  $\om$-rational functions. 
For instance an $\om$-automatic structure is defined via synchronous infinitary rational relations, see \cite{BlumensathGraedel04,KuskeLohrey}. On the other hand any synchronous infinitary rational relation is uniformizable by an $\om$-rational function, see \cite{ChoffrutGrigorieffG99}. Thus we can expect that our result will be useful in particular in the study of $\om$-automatic structures. 

\bibliographystyle{Alpha}
\bibliography{mabiblio,mabiblio-Set-theory}

\end{document}